\documentclass[draft,preprint,12pt,numbers,sort&compress]{elsarticle}
\usepackage{}
\usepackage{mathrsfs}
\usepackage{amssymb}
\usepackage{amsthm}
\usepackage{mathrsfs}
\usepackage[centertags]{amsmath}
\usepackage{amsfonts}

\usepackage{color}

\usepackage{ctex}
\usepackage{CJK}

\input amssym.def
\input amssym.tex

\date{}

\newtheorem{Theorem}{Theorem}[section]

\newtheorem{Lemma}{Lemma}[section]

\newcommand\R{\mbox{\bf R}}

\newcommand\Z{\mbox{\bf Z}}
\newcommand\z{\mbox{\bf z}}
\newcommand\SR{\mbox{\scriptsize\bf R}}

\newcommand{\definition}{{\lower .5ex
  \hbox{$\>\>\stackrel{\triangle}{=}\>\>$} }}
\newcommand\supp{\mathop{\rm supp}}

\begin{document}

\baselineskip=22pt
\thispagestyle{empty}

\begin{center}
{\Large \bf Probabilistic  pointwise convergence problem of some }\\[1ex]
{\Large \bf dispersive equations}\\[1ex]

{Wei Yan\footnote{Emails:   011133@htu.edu.cn}$^{a}$, \quad
Jinqiao Duan \footnote{Emails:  duan@iit.edu}$^{b}$,\quad
Yongsheng Li\footnote{Emails:
yshli@scut.edu.cn. }$^{c}$, \quad Meihua Yang\footnote{Corresponding author,Email: yangmeih@hust.edu.cn}$^{d*}$}\\[1ex]

{$^a$School of Mathematics and Information Science, Henan
Normal University,}\\
{Xinxiang, Henan 453007,   China}\\[1ex]

{$^b$Department of Applied Mathematics, Illinois Institute of Technology,}\\[1ex]
{Chicago, IL 60616, USA}\\[1ex]

{$^c$School of Mathematics, South  China  University of  Technology,}\\
{Guangzhou, Guangdong 510640,    China}\\[1ex]

{$^d$School of mathematics and statistics, Huazhong University of Science and Technology, }\\
{Wuhan, Hubei 430074,  China}\\[1ex]
\end{center}
\noindent{\bf Abstract.} In this paper, we investigate the almost surely pointwise convergence
 problem of free KdV equation, free wave equation, free elliptic and non-elliptic Schr\"odinger
  equation respectively. We firstly  establish some estimates
   related to the Wiener decomposition of frequency spaces  which is just Lemmas 2.1-2.6 in this paper.
    Secondly, by using  Lemmas 2.1-2.6, 3.1,  to establish the  probabilistic estimates
    of some random series which is just Lemmas 3.2-3.11 in this paper. Finally, combining  the
      density theorem in $H^{s}$ with Lemmas 3.2-3.11,  we obtained almost
       surely pointwise convergence of the solutions to
     corresponding equations with randomized initial data in $L^{2}$, which require much
      less regularity of the initial data than the rough data case. We also present the
       probabilistic density theorem.
 \medskip


\noindent {\bf Corresponding Author:}Meihua Yang

\noindent {\bf Email Address:}yangmeih@hust.edu.cn


\leftskip 0 true cm \rightskip 0 true cm

\newpage

\baselineskip=20pt

\bigskip
\bigskip

{\large\bf 1. Introduction}
\bigskip

\setcounter{Theorem}{0} \setcounter{Lemma}{0}\setcounter{Definition}{0}\setcounter{Proposition}{1}

\setcounter{section}{1}

In this paper,  we  investigate the pointwise convergence  problem of the free KdV equation in $\mathbf{R}$
\begin{eqnarray}\label{1.01}
\begin{cases}
u_t+\partial_{x}^{3}u=0, \ (x,t)\in\mathbf{R}\times\mathbf{R},\\
u(x,0)=f(x),     \   x\in\mathbf{R},
\end{cases}
\end{eqnarray}
free wave equation in $\mathbf{R}^n$, $n\geq 2$,
 \begin{eqnarray}\label{1.03}
\begin{cases}
u_{tt}+\Delta u=0, \ (x,t)\in\mathbf{R^{n}}\times\mathbf{R},\\
u(x,0)=f(x),  u_{t}(x,0)=0,   \   x\in\mathbf{R^{n}},
\end{cases}
\end{eqnarray}
and
free Schr\"odinger equation in $\mathbf{R}^n$, $n\geq 1$,
\begin{eqnarray}\label{1.04}
\begin{cases}
iu_t+\Delta_{\pm} u=0, \ (x,t)\in\mathbf{R^{n}}\times\mathbf{R},\\
u(x,0)=f(x),     \   x\in\mathbf{R^{n}}.
\end{cases}
\end{eqnarray}
Here $\Delta_{\pm}=\sum\limits_{j=1}^{n}\epsilon_{j}\partial_{x_{j}}^{2}, \epsilon_{j}=\pm 1.$
The formal solutions to the free KdV \eqref{1.01}, the free wave equation \eqref{1.03} and the
free Schr\"odinger equation \eqref{1.03} are given respectively by
\begin{equation}\label{S_1}
S_{1}(t)f(x)=(2\pi)^{-\frac{1}{2}}\int_{\SR}e^{ix\xi+it\xi^{3}}\mathscr{F}_{x}f(\xi)d\xi,
\end{equation}
\begin{equation}\label{S_2}
S_{2\pm}(t)f(x)=(2\pi)^{-\frac{n}{2}}\int_{\SR^{n}}
e^{ix\xi\pm it|\xi|}\mathscr{F}_{x}f(\xi)d\xi_{1}d\xi_{2}\cdot\cdot\cdot d\xi_{n},
\end{equation}
and
\begin{equation}\label{S_3}
S_{3}(t)f(x_{1},x_{2},\cdot\cdot\cdot,x_{n})=(2\pi)^{-\frac{n}{2}}
\int_{\SR^{n}}e^{ix\xi+it\left[\sum\limits_{j=1}^{n}\epsilon_{j}\xi_{j}^{2}\right]}
\mathscr{F}_{x}f(\xi)d\xi_{1}d\xi_{2}\cdot\cdot\cdot d\xi_{n},\epsilon_{j}=\pm 1,
\end{equation}
where
\begin{eqnarray*}
&&\mathscr{F}_xf(\xi)=(2\pi)^{-\frac{1}{2}}
\int_{\SR}e^{-ix\xi}f(x)dx,\\
&&\mathscr{F}_xf(\xi_{1},\xi_{2},\cdot\cdot\cdot,\xi_{n})=(2\pi)^{-\frac{n}{2}}
\int_{\SR^{n}}e^{-i\sum\limits_{j=1}^{n}x_{j}\xi_{j}}f(x)dx_{1}dx_{2}\cdot\cdot\cdot dx_{n}.
\end{eqnarray*}

The pointwise problem was originally studied by  Carleson \cite{Carleson}, who
showed pointwise convergence problem of the one dimensional Schr\"odinger equation in $H^{s}(\R)$, $s\geq 1/4$.
  The necessary condition and sufficient  condition for  the  pointwise convergence problem of the Schr\"odinger equation
attracts much attentions. For instance, Dahlberg and Kenig \cite{DK} showed that $s\geq \frac{1}{4}$ is
the necessary condition for the pointwise
convergence problem of the Schr\"odinger equation in any dimension.Dahlberg,  Kenig \cite{DK} and Kenig  et al. \cite{KPV1991,KPV1993}
  have showed the pointwise convergence problem
of  KdV equation in $H^{s}(\R)$ if and only if  $s\geq \frac{1}{4}.$ Bourgain \cite{Bourgain2016} presented counterexamples about  Schr\"odinger
 equation showing
 that convergence can
 fail if $s<\frac{n}{2(n+1)}$. Du et al. \cite{DGL} proved
 that the pointwise convergence problem
  of two dimensional
  Schr\"odinger
 equation in $H^{s}(\R)$ with $s>\frac{1}{3}.$ Du and Zhang \cite{DZ}
  proved  the
  pointwise convergence problem of $n$ dimensional Schr\"odinger
 equation in $H^{s}(\R)$ with $s>\frac{n}{2(n+1)},n\geq3.$ Thus, $\frac{n}{2(n+1)},n\geq2$ is optimal for
 the  pointwise convergence problem of the Schr\"odinger equation.
 Associated to the wave equation, Rogers and Villarroya \cite{RV} have proved that
$\frac{1}{2}\left[e^{it\sqrt{-\Delta}}+e^{-it\sqrt{-\Delta}}\right]
f\longrightarrow f$ with $f\in H^{s}(\R^{n})$
if and only if $s>{\rm max}\left\{n(\frac{1}{2}-\frac{1}{q}),
\frac{n+1}{4}-\frac{n-1}{2q},\frac{1}{2}\right\}(q\geq1)$.
 For the pointwise
 convergence problem of the  Schr\"odinger equation  in higher  dimension and
 other dispersive equations,
we also refer the readers to
\cite{Bourgain1992,Bourgain1995,Bourgain2013,CLV,Cowling,DG,DGZ,DK,GS,MVV,Lee,
LR2015,LR2017,MYZ2015,MZZ2015,KPV1991,KPV1993,RVV,Shao,S,Tao,TV,Vega}.

Recently, Compaan et al. \cite{CLS} applied randomized initial data to
 study pointwise convergence  of the Schr\"odinger flow, and then prove
 almost everywhere convergence
   with less regularity of the initial data. The method of  the suitably
     randomized initial data originated
from Lebowitz-Rose-Speer \cite{LRS} and Bourgain \cite{Bourgain1994, Bourgain1996}
 and Burq-Tzvetkov \cite{BTL,BTG}.
Many authors applied the method to study nonlinear dispersive
equations and hyperbolic equations in scaling super-critical regimes, for example,
see \cite{BOP,BOP2015,BOP2017,CG,CO,CZ,DC,DG,DLM,HM,KM,LM,NORS,NPS,NS,OOP,OP,P,ZF2011,ZF2012}.

In this paper, inspired by \cite{CLS, WH}, we mainly investigate the almost surely pointwise convergence problem
 of free KdV equation, free wave equation and elliptic and non-elliptic Schr\"odinger
  equation with randomized initial data in $L^2$, respectively. The main tools that we
  used are the density theorem and some estimates related to
   the Wiener decomposition of  the frequency spaces and Lemma 3.1. The crucial ingredients introduced in this paper
 are the  probabilistic estimates of some  random series which are just Lemmas 3.2-3.11 in this paper.

We give some notations before presenting our main results.
For $x\in \R^{n},$ we define $x^{\alpha}=\prod\limits_{j=1}^{n}x_{j}^{\alpha_{j}}$,
$\partial^{\beta}\phi
=\prod\limits_{j=1}^{n}(\partial/\partial_{x_{j}})^{\beta_{j}}\phi$,
 where $\alpha=\sum\limits_{j=1}^{n}\alpha_{j}, \beta=\sum\limits_{j=1}^{n}\beta_{j}.$
 For $\xi \in \R^{n},$ we have $|\xi|=\sqrt{\sum\limits_{j=1}^{2}\xi_{j}^{2}}.$
We define
\begin{eqnarray*}
&&D_{x}^{a}S_{1}(t)f=\frac{1}{\sqrt{2\pi}}\int_{\SR}e^{ix\xi+it\xi^{3}}
|\xi|^{a}\mathscr{F}_{x}f(\xi)d\xi,\\
&&D_{x}^{a}S_{2\pm}(t)f=\frac{1}{(2\pi)^{\frac{n}{2}}}\int_{\SR^{n}}
e^{ix\xi\pm it|\xi|}|\xi|^{a}\mathscr{F}_{x}f(\xi)d\xi,\\
&&D_{x}^{a}S_{3}(t)f=\frac{1}{(2\pi)^{\frac{n}{2}}}\int_{\SR^{n}}
e^{ix\xi+it(\sum_{j=1}^{n}\epsilon_{j}\partial_{x_{j}}^{2})}
|\xi|^{2a}\mathscr{F}_{x}f(\xi)d\xi,\\
&&D_{t}^{a}S_{1}(t)f=\frac{1}{\sqrt{2\pi}}\int_{\SR}e^{ix\xi+it\xi^{3}}
|\xi|^{3a}\mathscr{F}_{x}f(\xi)d\xi,\\
&&D_{t}^{a}S_{2\pm}(t)f=\frac{1}{(2\pi)^{\frac{n}{2}}}\int_{\SR^{n}}
e^{ix\xi\pm it|\xi|}|\xi|^{a}\mathscr{F}_{x}f(\xi)d\xi,\\
&&D_{t}^{a}S_{3}(t)f=\frac{1}{(2\pi)^{\frac{n}{2}}}\int_{\SR^{n}}
e^{ix\xi\pm it(\sum_{j=1}^{n}\epsilon_{j}\partial_{x_{j}}^{2})}|\sum_{j=1}^{n}\epsilon_{j}\partial_{x_{j}}^{2}|^{2a}\mathscr{F}_{x}f(\xi)d\xi.
\end{eqnarray*}
Now we introduce the randomization procedure for the initial data,
which can be seen in \cite{BOP,BOP2015,LM,ZF2012}.Let $B(0,1)$ be a unit ball centered in zero with radius equal to 1.
Let  $\psi\in C_{c}^{\infty}(\R^{n})$ be a real-valued,
even, non-negative bump function with $\supp \psi\subset B(0,1)$ such that
$\sum\limits_{k\in Z^{n}}\psi(\xi-k)=1$ for all $\xi \in \R^{n},$
which is known as  Wiener decomposition of  the frequency space.
For $s\in \R$ and  $f\in H^{s}(\R^{n})$.  For every $\xi \in Z^{n}$, we define the function $\psi(D-k)f:\R^n\rightarrow\mathbb{C}$ by
\begin{eqnarray}
(\psi(D-k)f)(x)=\mathscr{F}^{-1}\big(\psi(\xi-k)\mathscr{F}f\big)(x),x\in \R^n\label{1.010}.
\end{eqnarray}
If $f\in H^{s}$ for some $s\in \R$, then $P(D-k)f\in  H^{s}$ and
\begin{eqnarray}
f=\sum\limits_{k\in Z^{n}}P(D-k)f\label{1.011}
\end{eqnarray}
in $H^{s}$ with
\begin{eqnarray*}
\|f\|_{H^{s}}\sim \left[\sum\limits_{k\in Z^{n}}\|P_{k}f\|_{H^{s}}^{2}\right]^{\frac{1}{2}}= \left[\sum\limits_{k\in Z^{n}}\|P(D-k)f\|_{H^{s}}^{2}\right]^{\frac{1}{2}}.
\end{eqnarray*}
We will crucially exploit that these projections satisfy
 a unit-scale Bernstein inequality,
 namely that for all
$2 \leq p_1 \leq p_2 \leq\infty$
there exists $C\equiv C(p_1, p_2)>0$ such that for all $f\in L_{x}^2(
\R^n)$
and for all $k\in Z^n$
 \begin{eqnarray}\label{1.012}
\left\|\psi(D-k)f\right\|_{L_{x}^{p_2}(\SR^n)}\leq C\left\|\psi(D-k)f\right\|_{L_{x}^{p_1}(\SR^n)}\leq C\left\|\psi(D-k)f\right\|_{L_{x}^{2}(\SR^n)}.
\end{eqnarray}
Let $\{g_k\}_{k\in Z
^n}$ be a sequence of independent, zero-mean, complex-valued Gaussian random variables
 on a probability space $(\Omega,\mathcal{A}, \mathbb{P})$,
 where the real and imaginary parts of $g_k$ are independent and endowed with probability
distributions $\mu_k^1$ and $\mu_k^2$, respectively.
 Assume that there exists $c>0$ such that
 \begin{eqnarray}\label{1.013}
\Big|\int_{-\infty}^{+\infty}e^{\gamma x}d\mu_k^j(x)\Big|\leq e^{c\gamma^2},
\end{eqnarray}
for all $\gamma\in \R$, $k\in  \Z^n$, $j=1, 2$.
Thereafter for a given $f\in H^{s}(\R^n)$, $n\geq 1$, we define its
randomization by
\begin{eqnarray}
f^\omega:=\sum_{k\in \z^n}g_k(\omega)\psi(D-k)f.\label{1.014}
\end{eqnarray}
Lemma B.1 in \cite{BTL} showed that there is no smoothing upon
randomization in terms of differentiability. This randomization  improved the
 integrability of $f$, see Lemma 2.3 of \cite{BOP}. Such results for random Fourier series are known as Paley-Zygmund's theorem \cite{PZ}.
We define
\begin{eqnarray*}
\|f\|_{L_{\omega}^{p}(\Omega)}=\left[\int_{\Omega}|f(\omega)|^{p}dP(\omega)\right]^{\frac{1}{p}}.
\end{eqnarray*}
Obviously, $\|\|f^\omega\|_{H^{s}}\|_{L_{\omega}^{2}}=\|f\|_{H^{s}}$. We will restrict ourselves to a subset $\sum\subset \Omega$ with
$P\left(\sum\right)=1$ such that $f^{\omega}\in H^{s}$ for every $\omega \in \Omega.$

Then we show the main results of this paper as following:

 \begin{Theorem} \label{Theorem1}
Let  $f\in L^{2}(\R)$ and $f^{\omega}$ be a randomization of $f$ as
defined in (\ref{1.014}). Then, $\forall \alpha>0,$  we have
\begin{eqnarray}
\lim\limits_{t\longrightarrow0}{\rm P}\left(\omega\in \Omega: |S_{1}(t)f^{\omega}-f^{\omega}|>\alpha\right)=0. \label{1.016}
\end{eqnarray}

\end{Theorem}

\noindent {\bf Remark 1.} Dahlberg,  Kenig \cite{DK} and Kenig  et al. \cite{KPV1991,KPV1993}
  have showed the pointwise convergence problem
of  KdV equation in $H^{s}(\R)$ if and only if  $s\geq \frac{1}{4}.$ Obviously,
\begin{eqnarray}
\lim\limits_{\epsilon\longrightarrow 0}\alpha=\lim\limits_{\epsilon\longrightarrow 0}Ce\epsilon\left[{\rm In}\frac {C_{2}}{\epsilon}\right]^{\frac{1}{2}}=0\label{1.018}
\end{eqnarray}
and $\alpha=o(\epsilon^{\frac{1}{2}}).$
From \cite{DK,KPV1991,KPV1993} and Theorem 1.1,  we know that  the pointwise convergence problem of KdV equation
 with  random data requires less regularity of the initial data than
the pointwise convergence problem of KdV equation with  rough data.

\begin{Theorem} \label{Theorem2}
Let $f\in  L^{2}(\R^{n})$ and $f^{\omega}$ be a randomization of $f$ as
defined in (\ref{1.014}). Then, $\forall \alpha >0,$ we have
\begin{eqnarray}
\lim\limits_{t\longrightarrow 0}{\rm P}\left(\omega\in \Omega: \left|\frac{1}{2}\left[S_{2+}(t)
f^{\omega}(x)+S_{2-}
f^{\omega}(x)\right]-f^{\omega}(x)\right|>\alpha\right)\label{1.019}.
\end{eqnarray}

\end{Theorem}

\noindent{\bf Remark 2.} Rogers and Villarroya \cite{RV} have proved that
$\frac{1}{2}\left[e^{it\sqrt{-\Delta}}+e^{-it\sqrt{-\Delta}}\right]f\longrightarrow f$ with $f\in H^{s}(\R^{n})$
if and only if $s>{\rm max}\left\{n(\frac{1}{2}-\frac{1}{q}),\frac{n+1}{4}-\frac{n-1}{2q},\frac{1}{2}\right\}(q\geq1).$
From \cite{RV} and Theorem 1.2,  we know that  the pointwise convergence problem of wave equation
 with  random data requires less regularity of the initial data than
the pointwise convergence problem of wave equation with  rough data.

\begin{Theorem} \label{Theorem3}
Let $f\in  L^{2}(\R^{n})$ and  $f^{\omega}$ be a randomization of $f$ as
defined in (\ref{1.014}). Then, $\forall \alpha >0,$  we have
\begin{eqnarray}
\lim\limits_{t\longrightarrow 0}{\rm P}\left(\omega\in \Omega: \left|S_{3}(t)f^{\omega}-f^{\omega}\right|>\alpha\right)=0.\label{1.022}
\end{eqnarray}

\end{Theorem}

\noindent {\bf Remark 3.} Compaan et al. \cite{CLS} have proved the almost surely pointwise convergence problem in $H^{s}(s>0)$ for elliptic
Schr\"odinger equation with random data. Thus, our result improves the result of
\cite{CLS} to elliptic and non-elliptic Schr\"odinger equation. From \cite{Bourgain2016,D,DGL,DZ} and Theorem 1.3,
we know that  the pointwise convergence problem of elliptic Schr\"odinger equation
 with  random data requires less regularity of the initial data than
the pointwise convergence problem of elliptic Schr\"odinger equation  with  rough data.
Rogers et al. \cite{RVV}  showed that the solution to the two dimensional non-elliptic
Schr\"oodinger equation  converges to its initial datum $f$,
for all  $f\in H^{s}(\R^{2})$ if and only if $s \geq \frac{1}{2}.$ Thus, from \cite{RVV} and Theorem 1.3, we
know that  the pointwise convergence problem of two dimensional non-elliptic Schr\"odinger equation
 with  random data requires less regularity of the initial data than
the pointwise convergence problem of two dimensional non-elliptic Schr\"odinger equation  with  rough data.

\begin{Theorem} \label{Theorem4}(Probabilistic density Theorem)
 For $f\in L^{2}(\R^{n})$ and $\forall \epsilon>0,$  there exist a rapidly decreasing
 function $g$   and  $h\in L^{2}(\R)$ with $\|h\|_{L^{2}(\SR)}<\epsilon$ such that
 \begin{eqnarray*}
 \forall\omega\in \Omega_{\lambda M}:=\left\{\omega\in \Omega: \|h^{\omega}\|_{L^{2}}\leq \lambda\right\}\cap \left\{\omega\in \Omega: \left|x^{\alpha}\partial^{\beta}g^{\omega}\right|\leq M\right\},
 \end{eqnarray*}
 we have  $f^{\omega}=g^{\omega}+h^{\omega}$. Here,
 \begin{eqnarray*}
\lambda:=Ce\epsilon\left({\rm In}\frac{C_{1}}{\epsilon}\right)^{\frac{1}{2}},  M:=Ce\left[{\rm In}\frac{C_{1}}{\epsilon}\right]^{\frac{1}{2}}.
\end{eqnarray*}
and
\begin{eqnarray*}
\mathbb{P}\left(\Omega_{\lambda M}\right)\geq 1-2\epsilon.
\end{eqnarray*}

\end{Theorem}

Now, we present the outline of proof of Theorem 1.1 to explain the main idea of this paper
 since the Theorem 1.2, 1.3 can be proved similarly to Theorem 1.1.

 More precisely,  $f\in L^{2}$ and  since
rapidly decreasing  functions are dense in $L^2$(the density theorem which is just Lemma 2.2 in \cite{D}),
 we write $f=g+h$, where $g$ is a rapidly
 decreasing function and  $\|h\|_{L^{2}}<\epsilon$. Since $f^{\omega}=g^{\omega}+h^{\omega}$, then we get
\begin{eqnarray}
&&S_{1}(t)f^{\omega}-f^{\omega}=S_{1}(t)g^{\omega}-g^{\omega}
+S_{1}(t)h^{\omega}-h^{\omega}.\label{1.025}
\end{eqnarray}
Here, $f^{\omega}$ is defined as in (\ref{1.014}).

Then, when $|t|<\epsilon$,  $\alpha=Ce\epsilon\left[{\rm In}\frac {3C_{1}}
{\epsilon}\right]^{\frac{1}{2}}$,  we have
\begin{eqnarray*}
&&\mathbb{P}\left(\left\{\omega\in \Omega:
\left|S_{1}(t)f^{\omega}-f^{\omega}\right|>\alpha\right\}\right)\nonumber\\&&\leq
\mathbb{P}\left(\left\{\omega\in \Omega:
\left|S_{1}(t)g^{\omega}-g^{\omega}\right|>\frac{\alpha}{2}\right\}\right)
+\mathbb{P}\left(\left\{\omega\in \Omega:
\left|S_{1}(t)h^{\omega}\right|>\frac{\alpha}{4}\right\}\right)
\nonumber\\&&\qquad+\mathbb{P}\left(\left\{\omega\in \Omega: |h^{\omega}|>\frac{\alpha}{4}\right\}\right).
\end{eqnarray*}
Hence, we only need to deal with the right-hand side terms of the above inequality one by one.
 Note that $g$  is  a  rapidly  decreasing  function and $\|h\|_{L^2}<\epsilon$,
 and then combining the  probabilistic estimate Lemma 3.1
 with Strichartz estimates related to the uniform partition to the frequency spaces,
   we obtained the following estimates, the proofs are given in Lemma 3.2, Lemma 3.3 and Lemma 3.8, respectively.
 \begin{eqnarray}
\mathbb{P}\left(\left\{\omega\in \Omega:
\left|S_{1}(t)g^{\omega}-g^{\omega}\right|>\frac{\alpha}{2}\right\}\right)\leq
C_{1}e^{-\left(\frac{\alpha}{C|t|e}\right)^{2}},\label{1.026}
\end{eqnarray}
\begin{eqnarray}
\mathbb{P}\left(\left\{\omega\in \Omega:
\left|S_{1}(t)h^{\omega}\right|>\frac{\alpha}{4}\right\}\right)\leq C_{1}e^{-\left(\frac{\alpha}{Ce\|h\|_{H^{8\epsilon}}}\right)^{2}}\leq
 C_{1}e^{-\left(\frac{\alpha}{Ce\epsilon}\right)^{2}}, \label{1.027}
\end{eqnarray}
and
\begin{eqnarray}
\mathbb{P}\left(\left\{\omega\in \Omega: |h^{\omega}|>
\frac{\alpha}{4}\right\}\right)\leq C_{1}e^{-\left(\frac{\alpha}
{Ce\|h\|_{H^{8\epsilon}}}\right)^{2}}\leq
 C_{1}e^{-\left(\frac{\alpha}{Ce\epsilon}\right)^{2}}\label{1.028}.
\end{eqnarray}
Thus, when $|t|<\epsilon$,  $Ce\epsilon\left[{\rm In}
\frac {3C_{1}}{\epsilon}\right]^{\frac{1}{2}}\leq\alpha$,   we have
\begin{eqnarray}
&&\mathbb{P}\left(\left\{\omega\in \Omega:
\left|S_{1}(t)f^{\omega}-f^{\omega}\right|>\alpha\right\}\right)\nonumber\\&&\leq
\mathbb{P}\left(\left\{\omega\in \Omega:
\left|S_{1}(t)g^{\omega}-g^{\omega}\right|>\frac{\alpha}{2}\right\}\right)+\mathbb{P}\left(\left\{\omega\in \Omega:
\left|S_{1}(t)h^{\omega}\right|>\frac{\alpha}{4}\right\}\right)
\nonumber\\&&\qquad+\mathbb{P}\left(\left\{\omega\in \Omega: |h^{\omega}|>\frac{\alpha}{4}\right\}\right)\nonumber\\&&
\leq C_{1}e^{-\left(\frac{\alpha}{C|t|e}\right)^{2}}+2C_{1}e^{-\left(\frac{\alpha}{Ce\epsilon}\right)^{2}}\leq 3C_{1}e^{-\left(\frac{\alpha}{Ce\epsilon}\right)^{2}}\leq
\epsilon.\label{1.029}
\end{eqnarray}

The proof of the remainder of Theorem 1.1 can be seen in Lemma 3.11.

This completes the proof of Theorem 1.1.

\bigskip

\bigskip

\setcounter{section}{2}

\noindent{\large\bf 2. Preliminaries }

\setcounter{equation}{0}

\setcounter{Theorem}{0}

\setcounter{Lemma}{0}

\setcounter{section}{2}
In this section, we give some estimates related to the Wiener decomposition of the frequency spaces.

\begin{Lemma}\label{lem2.1}
For $f\in L^{2}(\R^{n})$, we have
\begin{eqnarray}
\left[\sum\limits_{k\in Z^{n}}|\psi(D-k)f|^{2}\right]^{\frac{1}{2}}\leq
\left\|f\right\|_{L^{2}(\SR^{n})}\label{2.01}.
\end{eqnarray}
\end{Lemma}
\noindent{\bf Proof.} To obtain (\ref{2.01}), it suffices to prove
\begin{eqnarray}
\sum\limits_{k\in Z^{n}}|\psi(D-k)f|^{2}\leq
\left\|f\right\|_{L^{2}(\SR^{n})}^{2}\label{2.02}.
\end{eqnarray}
By using the Cauchy-Schwarz inequality with respect to $\xi$, since $\supp\psi\subset B(0,1)$ we have
\begin{eqnarray}
&&\sum\limits_{k\in Z^{n}}|\psi(D-k)f|^{2}=\frac{1}{(2\pi)^{\frac{n}{2}}}\sum\limits_{k\in Z^{n}}\left|\int_{\SR^{n}}e^{i\sum\limits_{j=1}^{n}x_{j}\xi_{j}}\psi(\xi-k)\mathscr{F}_{x}f(\xi)d\xi\right|^{2}\nonumber\\
&&=\frac{1}{(2\pi)^{\frac{n}{2}}}\sum\limits_{k\in Z^{n}}\left|\int_{|\xi-k|\leq1}e^{i\sum\limits_{j=1}^{n}x_{j}\xi_{j}}\psi(\xi-k)\mathscr{F}_{x}f(\xi)d\xi\right|^{2}\nonumber\\
&&\leq\left[\sum\limits_{k\in Z}\int_{|\xi-k|\leq1}|\psi(\xi-k)\mathscr{F}_{x}f(\xi)|^{2}d\xi\right]\left[\int_{|\xi-k|\leq1}d\xi\right]\nonumber\\
&&\leq \left[\sum\limits_{k\in Z}\int_{|\xi-k|\leq1}|\psi(\xi-k)\mathscr{F}_{x}f(\xi)|^{2}d\xi\right]\nonumber\\
&&=\sum\limits_{k\in Z^{n}}\left\|\psi(\xi-k)\mathscr{F}_{x}f(\xi)\right\|_{L^{2}}^{2}.\label{2.03}
\end{eqnarray}
From
\begin{eqnarray}
\mathscr{F}_{x}f(\xi)=\sum\limits_{k\in Z^{n}}\psi(\xi-k)\mathscr{F}_{x}f(\xi),\label{2.04}
\end{eqnarray}
by using $\supp\psi\subset B(0,1)$,  we have
\begin{eqnarray}
&&\|\mathscr{F}_{x}f(\xi)\|_{L^{2}}^{2}=
\sum\limits_{k\in Z^{n}}\sum\limits_{l\in Z^{n}}\int_{\SR}\left[\psi(\xi-k)\mathscr{F}_{x}f(\xi)\right]\left[\psi(\xi-l)\overline{\mathscr{F}_{x}f}(\xi)\right]d\xi\nonumber\\&&=
\sum\limits_{k\in Z^{n}}\int_{\SR}\left|\psi(\xi-k)\mathscr{F}_{x}f(\xi)\right|^{2}d\xi\label{2.05}.
\end{eqnarray}
Combining (\ref{2.03}) with (\ref{2.05}), we derive (\ref{2.02}).

This completes the proof of Lemma \ref{lem2.1}. $\hfill\Box$

\begin{Lemma}\label{lem2.2}
For $f\in L^{2}(\R)$, we have
\begin{eqnarray}
\left[\sum\limits_{k\in Z}|\psi(D-k)S_{1}(t)f|^{2}\right]^{\frac{1}{2}}\leq
\left\|f\right\|_{L^{2}(\SR)}\label{2.06}.
\end{eqnarray}
\end{Lemma}
\noindent{\bf Proof.} To obtain (\ref{2.06}), it suffices to prove
\begin{eqnarray}
\sum\limits_{k\in Z}|\psi(D-k)S_{1}(t)f|^{2}\leq
\left\|f\right\|_{L^{2}(\SR)}^{2}\label{2.07}.
\end{eqnarray}
By using the Cauchy-Schwarz inequality with respect to $\xi$,  since $\supp \psi \subset B(0,1)$, we have
\begin{eqnarray}
&&\sum\limits_{k\in Z}|\psi(D-k)S_{1}(t)f|^{2}=\frac{1}{(2\pi)^{\frac{n}{2}}}\sum\limits_{k\in Z^{n}}\left|\int_{\SR}e^{ix\xi}e^{it\xi^{3}}\psi(\xi-k)\mathscr{F}_{x}f(\xi)d\xi\right|^{2}\nonumber\\
&&=\frac{1}{(2\pi)^{\frac{n}{2}}}\sum\limits_{k\in Z^{n}}\left|\int_{|\xi-k|\leq1}e^{ix\xi}e^{it\xi^{3}}\psi(\xi-k)\mathscr{F}_{x}f(\xi)d\xi\right|^{2}\nonumber\\
&&\leq \left[\sum\limits_{k\in Z}\int_{|\xi-k|\leq1}|\psi(\xi-k)\mathscr{F}_{x}f(\xi)|^{2}d\xi\right]\left[\int_{|\xi-k|\leq1}d\xi\right]\nonumber\\
&&\leq \left[\sum\limits_{k\in Z}\int_{|\xi-k|\leq1}|\psi(\xi-k)\mathscr{F}_{x}f(\xi)|^{2}d\xi\right]\nonumber\\
&&=\sum\limits_{k\in Z}\left\|\psi(\xi-k)\mathscr{F}_{x}f(\xi)\right\|_{L^{2}}^{2}.\label{2.08}
\end{eqnarray}
From
\begin{eqnarray}
\mathscr{F}_{x}f(\xi)=\sum\limits_{k\in Z}\psi(\xi-k)\mathscr{F}_{x}f(\xi),\label{2.09}
\end{eqnarray}
by using $\supp\psi\subset B(0,1)$,  we have
\begin{eqnarray}
&&\|\mathscr{F}_{x}f(\xi)\|_{L^{2}}^{2}=
\sum\limits_{k\in Z}\sum\limits_{l\in Z}\int_{\SR}\left[\psi(\xi-k)\mathscr{F}_{x}f(\xi)\right]\left[\psi(\xi-l)\overline{\mathscr{F}_{x}f}(\xi)\right]d\xi\nonumber\\&&=
\sum\limits_{k\in Z}\int_{\SR}\left|\psi(\xi-k)\mathscr{F}_{x}f(\xi)\right|^{2}d\xi\label{2.010}.
\end{eqnarray}
Combining (\ref{2.08}) with (\ref{2.010}), we derive (\ref{2.07}).

This completes the proof of Lemma \ref{lem2.2}. $\hfill\Box$

\begin{Lemma}\label{lem2.3}
For $f\in L^{2}(\R^{n})$, we have
\begin{eqnarray}
\left[\sum\limits_{k\in Z^{n}}|\psi(D-k)S_{2}(t)f|^{2}\right]^{\frac{1}{2}}\leq
\left\|f\right\|_{L^{2}(\SR^{n})}\label{2.011}.
\end{eqnarray}
\end{Lemma}

Lemma 2.3 can be proved similarly to Lemma 2.2.

\begin{Lemma}\label{lem2.4}
For $f\in L^{2}(\R^{n})$, we have
\begin{eqnarray}
\left[\sum\limits_{k\in Z^{n}}|\psi(D-k)S_{3}(t)f|^{2}\right]^{\frac{1}{2}}\leq
\left\|f\right\|_{L^{2}(\SR^{n})}\label{2.012}.
\end{eqnarray}
\end{Lemma}

Lemma 2.4 can be proved similarly to Lemma 2.2.
\begin{Lemma}\label{lem2.5}
Let $g$  be a rapidly decreasing function and we denote by $\psi^{(\beta)}$ the $\beta$ order derivative of $\psi$,  we have
\begin{eqnarray}
\sum_{|k|\geq3}\int_{\SR}\left|\xi^{\alpha}\mathscr{F}_{x}g(\xi)\psi^{(\beta)}(\xi-k)\right|^{2}d\xi\leq C.\label{2.013}
\end{eqnarray}
\end{Lemma}
\noindent{\bf Proof.} Since $\supp\psi\subset B(0,1)$, we have $\supp\psi^{(\beta)}\subset B(0,1)$.  Let $\xi-k=\eta$, then, $\xi=k+\eta$, since
 $g$  is a rapidly decreasing function,  we have
\begin{eqnarray}
&&\sum_{|k|\geq3}\int_{\SR}\left|\xi^{\alpha}\mathscr{F}_{x}g(\xi)\psi^{(\beta)}(\xi-k)\right|^{2}d\xi\nonumber\\&&=
\sum_{|k|\geq3}\int_{\SR}\left|(\eta+k)^{\alpha}\mathscr{F}_{x}g(\eta+k)\psi^{(\beta)}(\eta)\right|^{2}d\eta
\nonumber\\&&=\sum_{|k|\geq3}\int_{|\eta|\leq1}\left|(\eta+k)^{\alpha}\mathscr{F}_{x}g(\eta+k)\psi^{(\beta)}(\eta)\right|^{2}d\eta\nonumber\\
&&\leq \sum_{|k|\geq3}\int_{|\eta|\leq1}\frac{1}{1+\left|\eta+k\right|^{2}}d\eta\leq \sum_{|k|\geq3}\frac{C}{k^{2}}\leq C
.\label{2.014}
\end{eqnarray}

This completes the proof of Lemma 2.5.

\begin{Lemma}\label{lem2.6}
Let $g$  be a rapidly decreasing function,  we have
\begin{eqnarray}
\sum_{|k|\geq3}\int_{\SR^{n}}\left|\xi^{\alpha}\mathscr{F}_{x}g(\xi)\partial^{\beta}\psi(\xi-k)\right|^{2}d\xi\leq C.\label{2.015}
\end{eqnarray}
\end{Lemma}

Lemma 2.6 can be proved similarly to Lemma 2.5.

\bigskip
\bigskip

\noindent{\large\bf 3. Probabilistic estimates of some random series}

\setcounter{equation}{0}

\setcounter{Theorem}{0}

\setcounter{Lemma}{0}
\setcounter{section}{3}
In this section, we establish probabilistic estimates of some random series. More precisely,
we apply Lemmas 2.1-2.6 and 3.1 to establish Lemmas 3.2-3.11 which play  crucial role in establishing
Theorems 1.1-1.3. In particular,  we  apply Lemma 3.1 to establish Lemmas 3.9, 3.10 which can be used to establish  Lemma 3.11 which  is  called as the probabilistic
density theorem.

\begin{Lemma}\label{lem3.1}
Assume (\ref{1.013}). Then, there exists $C>0$ such that
\begin{eqnarray}
\left\|\sum_{k\in\z^n}g_k(\omega)c_k\right\|_{L_{\omega}^p(\Omega)}
\leq C\sqrt{p}\left\|c_k\right\|_{l^{2}(\z^n)}.\label{3.01}
\end{eqnarray}
for all $p\geq2 $
and $\{c_k\}\in l^{2}(\Z^n)$.
\end{Lemma}

For the proof of Lemma 3.1, we refer the readers to Lemma 3.1 of \cite{BTL}.

\begin{Lemma}\label{lem3.2}
 $\forall \alpha>0$. Let $g$ be  a  rapidly
decreasing  function
and we denote by $g^{\omega}$ the randomization of $g$ as
defined in (\ref{1.014}).
 Then, there exist $C>0,C_{1}>0$ such that
\begin{eqnarray}
\mathbb{P}\left(\Omega_{1}^{c}\right)\leq C_{1}e^{-\left(\frac{\alpha}{C|t|e}\right)^{2}}.\label{3.02}
\end{eqnarray}
where
\begin{eqnarray*}
\Omega_{1}^{c}=\left\{\omega\in \Omega:
\left|S_{1}(t)g^{\omega}-g^{\omega}\right|>\alpha\right\}.
\end{eqnarray*}
\end{Lemma}
\noindent{\bf Proof.}By  using Lemma 3.1 and  the Cauchy-Schwartz inequality,
 since $g$  is  a  rapidly  decreasing  function and $|e^{it\xi^{3}}-1|\leq |t\xi^{3}|,$
we have
\begin{eqnarray}
&&\left\|S_{1}(t)g^{\omega}-g^{\omega}\right\|_{L_{\omega}^{p}(\Omega)}\leq C\sqrt{p}\left[\sum_{k}\left|\int_{\SR}(e^{it\xi^{3}}-1)e^{ix \xi}\psi(\xi-k)\mathscr{F}g(\xi)d\xi\right|^{2}\right]^{\frac{1}{2}}\nonumber\\
&&\leq C|t|\sqrt{p}\left[\sum_{k}\int_{|\xi-k|\leq1}|\xi|^{3}\left|\psi(\xi-k)
\mathscr{F}g(\xi)\right|^{2}d\xi\right]^{\frac{1}{2}}\nonumber\\
&&\leq C|t|\sqrt{p}\left[\sum_{k}\int_{|\xi-k|\leq1}|\xi|^{6}
\left[\psi(\xi-k)\mathscr{F}g(\xi)\right]^{2}d\xi\right]^{\frac{1}{2}}\left[\int_{|\xi-k|\leq1}d\xi\right]^{\frac{1}{2}}\nonumber\\
&&\leq C|t|\sqrt{p}\left[\sum_{k}\int_{|\xi-k|\leq1}|\xi|^{6}
\left[\psi(\xi-k)\mathscr{F}g(\xi)\right]^{2}d\xi\right]^{\frac{1}{2}}\nonumber\\
&&=C|t|\sqrt{p}\left[\sum_{k}\left\|\psi(D-k)g\right\|_{H^{3}}^{2}\right]^{\frac{1}{2}}\nonumber\\
&&=C|t|\sqrt{p}\|g\|_{H^{3}}\leq C\sqrt{p}|t|.\label{3.03}
\end{eqnarray}
Thus, by using Chebyshev inequality, from (\ref{3.03}), we have
\begin{eqnarray}
&&\mathbb{P}\left(\Omega_{1}^{c}\right)\leq \int_{\Omega_{1}^{c}}\left[\frac{\left|S_{1}(t)g^{\omega}-g^{\omega}
\right|}{\alpha}\right]^{p}d\mathbb{P}(\omega)
\leq \left(\frac{C\sqrt{p}\|t|}{\alpha}\right)^{p}\label{3.04}.
\end{eqnarray}
Take
\begin{eqnarray}
p=\left(\frac{\alpha}{Ce|t|}\right)^{2}\label{3.05}.
\end{eqnarray}
If $p\geq2$, from (\ref{3.04}),
then we have
\begin{eqnarray}
&&\mathbb{P}\left(\Omega_{1}^{c}\right)\leq e^{-p}=e^{-\left(\frac{\alpha}{Ce|t|}\right)^{2}}\label{3.06}.
\end{eqnarray}
If $p\leq2$, from (\ref{3.04}),
 we have
\begin{eqnarray}
\mathbb{P}(\Omega_{1}^{c})\leq 1\leq e^{2}e^{-2}\leq C_{1}e^{-\left(\frac{\alpha}
{Ce|t|}\right)^{2}}.\label{3.07}
\end{eqnarray}
Here $C_{1}=e^{2}.$
Thus, from (\ref{3.06}), (\ref{3.07}), we have
\begin{eqnarray}
\mathbb{P}(\Omega_{1}^{c})\leq  C_{1} e^{
\left[-\left(\frac{\alpha}
{C|t|e}\right)^{2}\right]}.\label{3.08}
\end{eqnarray}

This completes the proof of Lemma 3.2.$\hfill\Box$

\begin{Lemma}\label{lem3.3}

 Let $h\in L^{2}(\R)$ and we denote by $h^{\omega}$
 the randomization of $h$ as defined in (\ref{1.014}).  Then, $\forall \alpha>0$
 there exist $C>0, C_{1}>0$ such that
\begin{eqnarray}
\mathbb{P}\left(\Omega_{2}^{c}\right)\leq C_{1}e^{-\left(\frac{\alpha}{Ce\|h\|_{L^{2}}}\right)^{2}}.\label{3.09}
\end{eqnarray}
where
$
\Omega_{2}^{c}=\left\{\omega\in \Omega:
\left|S_{1}(t)h^{\omega}\right|>\alpha\right\}.
$
\end{Lemma}
\noindent {\bf Proof.}
By using   Lemmas 3.1, 2.2,  we have
\begin{eqnarray}
&&\hspace{-1cm}\left\|S_{1}(t)h^{\omega}
\right\|_{L_{\omega}^{p}(\Omega)}
=\left\|\sum\limits_{k\in Z^{n}}g_{k}(\omega)\psi(D-k)S_{1}(t)h
\right\|_{L_{\omega}^{p}(\Omega)}
\nonumber\\
&&\hspace{-1cm}\leq C\sqrt{p}\left(\sum_{k\in Z^{n}}\left|\psi(D-k)S_{1}(t)h
\right|^{2}\right)^{\frac{1}{2}}
\leq C\sqrt{p}\|h\|_{L^{2}} \label{3.010}.
\end{eqnarray}
Thus, by using Chebyshev inequality, we have
\begin{eqnarray}
&&\mathbb{P}\left(\Omega_{2}^{c}\right)\leq \int_{\Omega_{2}^{c}}\left[\frac{\left|S_{1}(t)h^{\omega}
\right|}{\alpha}\right]^{p}d\mathbb{P}(\omega)
\leq \left(\frac{C\sqrt{p}\|h\|_{L^{2}}}{\alpha}\right)^{p}\label{3.011}.
\end{eqnarray}
By using a proof similar to (\ref{3.08}), we obtain (\ref{3.09}).

This completes the proof of Lemma 3.3.$\hfill\Box$

\bigskip
\bigskip

\begin{Lemma}\label{lem3.4}
  Let $g$ be  a  rapidly  decreasing  function
and we denote by $g^{\omega}$ the randomization of $g$ as defined in (\ref{1.014}).
 Then, $\forall\alpha>0$, there exist $C>0,C_{1}>0$ such that
\begin{eqnarray}
\mathbb{P}\left(\Omega_{3}^{c}\right)\leq C_{1}e^{-\left(\frac{\alpha}{C|t|e}\right)^{2}}.\label{3.012}
\end{eqnarray}
where
$
\Omega_{3}^{c}=\left\{\omega\in \Omega:
\left|\frac{1}{2}\left[S_{2+}(t)+S_{2-}(t)\right]g^{\omega}-g^{\omega}\right|>\alpha\right\}.
$
\end{Lemma}
\noindent{\bf Proof.} By  using Lemma 3.1 and  the Cauchy-Schwartz inequality with respect to $\xi$,
since   $g$ is is a  rapidly  decreasing  function and
$\left|\frac{1}{2}\left[e^{ it|\xi|}+e^{-it|\xi|}\right]-1\right)|\leq |t||\xi|,$ we have
\begin{eqnarray}
&&\left\|\frac{1}{2}\left[S_{2+}(t)+S_{2-}(t)\right]g^{\omega}-g^{\omega}\right\|_{L_{\omega}^{p}(\Omega)}\nonumber\\&&\leq C\sqrt{p}\left[\sum_{k}\left|\int_{\SR^{n}}\left(\frac{1}{2}\left[e^{ it|\xi|}+e^{-it|\xi|}\right]-1\right)e^{ix \xi}\psi(\xi-k)\mathscr{F}g(\xi)d\xi\right|^{2}\right]^{\frac{1}{2}}\nonumber\\
&&=\sqrt{p}\left[\sum_{k}\left|\int_{|\xi-k|\leq1}\left(\frac{1}{2}\left[e^{ it|\xi|}+e^{-it|\xi|}\right]-1\right)e^{ix \xi}\psi(\xi-k)\mathscr{F}g(\xi)d\xi\right|^{2}\right]^{\frac{1}{2}}\nonumber\\
&&\leq C|t|\sqrt{p}\left[\sum_{k}\left|\int_{|\xi-k|\leq1}\left|\xi\psi(\xi-k)
\mathscr{F}g(\xi)\right|d\xi\right|^{2}\right]^{\frac{1}{2}}\nonumber\\
&&\leq C|t|\sqrt{p}\left[\sum_{k}\int_{|\xi-k|\leq1}\left|\xi\psi(\xi-k)
\mathscr{F}g(\xi)\right|^{2}d\xi\right]^{\frac{1}{2}}\left[\int_{|\xi-k|\leq1}d\xi\right]^{\frac{1}{2}}\nonumber\\
&&\leq C|t|\sqrt{p}\left[\sum_{k}\int_{\SR^{n}}\left|\xi\psi(\xi-k)
\mathscr{F}g(\xi)\right|^{2}d\xi\right]^{\frac{1}{2}}\nonumber\\
&&=C|t|\sqrt{p}\left[\sum_{k}\left\|\psi(D-k)\right\|_{H^{1}}^{2}\right]^{\frac{1}{2}}\nonumber\\
&&=C|t|\sqrt{p}\left\|g\right\|_{H^{1}}\leq C|t|\sqrt{p}.\label{3.013}
\end{eqnarray}
Thus, from (\ref{3.013}), by using Chebyshev inequality, from (\ref{3.013}), we have
\begin{eqnarray}
\mathbb{P}(\Omega_{3}^{c})\leq \frac{\left\|\frac{1}{2}\left[S_{2+}(t)+S_{2-}(t)\right]g^{\omega}-g^{\omega}\right\|_{L_{\omega}^{p}(\Omega)}^{p}}{\alpha^{p}}\leq
\frac{(C|t|\sqrt{p})^{p}}{\alpha^{p}}.\label{3.014}
\end{eqnarray}
By using a  proof similar to (\ref{3.08}), from  (\ref{3.014}), we have
\begin{eqnarray}
\mathbb{P}(\Omega_{3}^{c})\leq  C_{1} {\rm exp}\left[-\left(\frac{\alpha}
{Ce|t|}\right)^{2}\right].\label{3.015}
\end{eqnarray}

This completes the proof of Lemma 3.4.$\hfill\Box$

\begin{Lemma}\label{lem3.5}
  Let $h\in L^{2}(\R^{n})$
 and we denote by $h^{\omega}$ the randomization of $h$ as defined in
 (\ref{1.014}).  Then,  $\forall \alpha>0$,   there exist $C>0$ and $C_{1}>0$ such that
\begin{eqnarray}
\mathbb{P}\left(\Omega_{4}^{c}\right)\leq C_{1}e^{-\left(\frac{\alpha}
{Ce\|h\|_{L^{2}}}\right)^{2}},\label{3.016}
\end{eqnarray}
where
\begin{eqnarray}
\Omega_{4}^{c}=\left\{\omega\in \Omega: \left|S_{2\pm}(t)h^{\omega}
\right|>\alpha\right\}.\label{3.017}
\end{eqnarray}
\end{Lemma}
\noindent {\bf Proof.} By using   Lemmas 3.1, 2.2,  we have
\begin{eqnarray}
&&\hspace{-1cm}\left\|S_{1}(t)h^{\omega}
\right\|_{L_{\omega}^{p}(\Omega)}
=\left\|\sum\limits_{k\in Z^{n}}g_{k}(\omega)\psi(D-k)S_{2\pm}(t)h
\right\|_{L_{\omega}^{p}(\Omega)}
\nonumber\\
&&\hspace{-1cm}\leq C\sqrt{p}\left(\sum_{k\in Z^{n}}\left|\psi(D-k)S_{2\pm}(t)h
\right|^{2}\right)^{\frac{1}{2}}
\leq C\sqrt{p}\|h\|_{L^{2}} \label{3.018}.
\end{eqnarray}
Thus, by using Chebyshev inequality, from (\ref{3.018}), we have
\begin{eqnarray}
&&\mathbb{P}\left(\Omega_{4}^{c}\right)\leq \int_{\Omega_{4}^{c}}\left[\frac{\left|S_{2\pm}(t)h^{\omega}
\right|}{\alpha}\right]^{p}d\mathbb{P}(\omega)
\leq \left(\frac{C\sqrt{p}\|h\|_{L^{2}}}{\alpha}\right)^{p}\label{3.019}.
\end{eqnarray}
By using a  proof similar to (\ref{3.08}), from  (\ref{3.019}), we have
\begin{eqnarray}
\mathbb{P}(\Omega_{4}^{c})\leq  C_{1} {\rm exp}\left[-\left(\frac{\alpha}
{Ce|t|}\right)^{2}\right].\label{3.020}
\end{eqnarray}

This completes the proof of Lemma 3.5.$\hfill\Box$

\begin{Lemma}\label{lem3.6}
  Let $g$ be   a  rapidly
decreasing  function
and we denote by $g^{\omega}$ the randomization of $g$ as
defined in (\ref{1.014}).
 Then, $\forall \alpha>0$,  there exist $C>0,C_{1}>0$ such that
\begin{eqnarray}
\mathbb{P}\left(\Omega_{5}^{c}\right)\leq C_{1}e^{-\left(\frac{\alpha}{C|t|e}\right)^{2}}.\label{3.021}
\end{eqnarray}
where
\begin{eqnarray*}
\Omega_{5}^{c}=\left\{\omega\in \Omega:
\left|S_{3}(t)g^{\omega}-g^{\omega}\right|>\alpha\right\}.
\end{eqnarray*}
\end{Lemma}
\noindent{\bf Proof.} By using Lemma 3.1 and  the Cauchy-Schwartz inequality with respect to $\xi$,
 since $g$  is  a  rapidly  decreasing  function and  $\left|e^{-it[\sum_{j=1}^{n}\epsilon_{j}\xi_{j}^{2}]}-1\right|\leq t|\xi|^{2}$. we have
\begin{eqnarray}
&&\left\|S_{3}(t)g^{\omega}-g^{\omega}\right\|_{L_{\omega}^{p}(\Omega)}\leq C\sqrt{p}\left[\sum_{k}\left|\int_{\SR^{n}}(e^{-it[\sum_{j=1}^{n}\epsilon_{j}\xi_{j}^{2}]}-1)e^{ix \xi}\psi(\xi-k)\mathscr{F}g(\xi)d\xi\right|^{2}\right]^{\frac{1}{2}}\nonumber\\
&&=C\sqrt{p}\left[\sum_{k}\left|\int_{|\xi-k|\leq1}(e^{-it[\sum_{j=1}^{n}\epsilon_{j}\xi_{j}^{2}]}-1)e^{ix \xi}\psi(\xi-k)\mathscr{F}g(\xi)d\xi\right|^{2}\right]^{\frac{1}{2}}\nonumber\\
&&\leq C|t|\sqrt{p}\left[\sum_{k}\int_{|\xi-k|\leq1}|\xi|^{2}\left|\psi(\xi-k)
\mathscr{F}g(\xi)\right|^{2}d\xi\right]^{\frac{1}{2}}\left[\int_{|\xi-k|\leq1}d\xi\right]^{\frac{1}{2}}\nonumber\\
&&= C|t|\sqrt{p}\left[\sum_{k}\int_{\SR^{n}}|\xi|^{2}|\psi(\xi-k)\mathscr{F}g(\xi)|^{2}d\xi
\right]^{\frac{1}{2}}\nonumber\\
&&\leq C|t|\sqrt{p}\left[\sum_{k}\left\|\psi(D-k)g\right\|_{H^{1}}^{2}\right]^{\frac{1}{2}}\nonumber\\
&&=C|t|\sqrt{p}\|g\|_{H^{1}}\leq C\sqrt{p}|t|.\label{3.022}
\end{eqnarray}
From (\ref{3.020}), by using  Chebyshev inequality, from (\ref{3.022}), we have
\begin{eqnarray}
\mathbb{P}(\Omega_{5}^{c})\leq \frac{(C\sqrt{p}|t|)^{p}}{\alpha^{p}}.\label{3.023}
\end{eqnarray}
Thus, by using a proof similar to (\ref{3.08}), from (\ref{3.023}),  we have
\begin{eqnarray}
\mathbb{P}(\Omega_{5}^{c})\leq  C_{1} {\rm exp}
\left[-\left(\frac{\alpha}
{C|t|e}\right)^{2}\right].\label{3.024}
\end{eqnarray}

This completes the proof of Lemma 3.6.$\hfill\Box$

 \begin{Lemma}\label{lem3.7}
 Let $h\in L^{2}(\R^{n})$
 and we denote by $h^{\omega}$ the randomization of $h$ as defined in
(\ref{1.014}).  Then,  $\forall \alpha>0$, there exist $C>0$ and $C_{1}>0$ such that
\begin{eqnarray}
\mathbb{P}\left(\Omega_{6}^{c}\right)\leq C_{1}^{-\left(\frac{\alpha}
{Ce\|h\|_{L^{2}}}\right)^{2}}.\label{3.025}
\end{eqnarray}
where
\begin{eqnarray}
\Omega_{6}^{c}=\left\{\omega\in \Omega: |S_{3}(t)h^{\omega}|>\alpha\right\}.\label{3.026}
\end{eqnarray}
\end{Lemma}
\noindent {\bf Proof.}By using   Lemmas 3.1, 2.2,  we have
\begin{eqnarray}
&&\hspace{-1cm}\left\|S_{3}(t)h^{\omega}
\right\|_{L_{\omega}^{p}(\Omega)}
=\left\|\sum\limits_{k\in Z^{n}}g_{k}(\omega)\psi(D-k)S_{2\pm}(t)h
\right\|_{L_{\omega}^{p}(\Omega)}
\nonumber\\
&&\hspace{-1cm}\leq C\sqrt{p}\left(\sum_{k\in Z^{n}}\left|\psi(D-k)S_{2\pm}(t)h
\right|^{2}\right)^{\frac{1}{2}}
\leq C\sqrt{p}\|h\|_{L^{2}} \label{3.027}.
\end{eqnarray}
Thus, by using Chebyshev inequality, from  (\ref{3.027}), we have
\begin{eqnarray}
&&\mathbb{P}\left(\Omega_{6}^{c}\right)\leq \int_{\Omega_{6}^{c}}\left[\frac{\left|S_{3}(t)h^{\omega}
\right|}{\alpha}\right]^{p}d\mathbb{P}(\omega)
\leq \left(\frac{C\sqrt{p}\|h\|_{L^{2}}}{\alpha}\right)^{p}\label{3.028}.
\end{eqnarray}
By using a  proof similar to (\ref{3.08}), from  (\ref{3.028}), we have
\begin{eqnarray}
\mathbb{P}(\Omega_{6}^{c})\leq  C_{1} {\rm exp}\left[-\left(\frac{\alpha}
{Ce|t|}\right)^{2}\right].\label{3.029}
\end{eqnarray}

This completes the proof of Lemma 3.7.$\hfill\Box$

\bigskip

\begin{Lemma}\label{lem3.8}
 Let $h\in L^{2}(\R^{n})$
 and we denote by $h^{\omega}$ the randomization of $h$ as defined in
(\ref{1.014}).  Then, $\forall \alpha>0$, there exist $C>0$ and $C_{1}>0$ such that
\begin{eqnarray}
\mathbb{P}\left(\Omega_{7}^{c}\right)\leq C_{1}
e^{-\left(\frac{\alpha}{Ce\|h\|_{L^{2}}}\right)^{2}}\label{3.030}.
\end{eqnarray}
where
\begin{eqnarray}
\Omega_{7}^{c}=\left\{\omega\in \Omega: |h^{\omega}|>\alpha\right\}\label{3.031}.
\end{eqnarray}
\end{Lemma}
\noindent{\bf Proof.}By using   Lemmas 3.1, 2.2,  we have
\begin{eqnarray}
&&\hspace{-1cm}\left\|h^{\omega}
\right\|_{L_{\omega}^{p}(\Omega)}
=\left\|\sum\limits_{k\in Z^{n}}g_{k}(\omega)\psi(D-k)h
\right\|_{L_{\omega}^{p}(\Omega)}
\nonumber\\
&&\hspace{-1cm}\leq C\sqrt{p}\left(\sum_{k\in Z^{n}}\left|\psi(D-k)h
\right|^{2}\right)^{\frac{1}{2}}
\leq C\sqrt{p}\|h\|_{L^{2}} \label{3.032}.
\end{eqnarray}
Thus, by using Chebyshev inequality, from (\ref{3.032}),  we have
\begin{eqnarray}
&&\mathbb{P}\left(\Omega_{7}^{c}\right)\leq \int_{\Omega_{7}^{c}}\left[\frac{\left|h^{\omega}
\right|}{\alpha}\right]^{p}d\mathbb{P}(\omega)
\leq \left(\frac{C\sqrt{p}\|h\|_{L^{2}}}{\alpha}\right)^{p}\label{3.033}.
\end{eqnarray}
By using a proof similar to (\ref{3.08}), we obtain (\ref{3.030}).

This completes the proof of Lemma 3.8.$\hfill\Box$

\begin{Lemma}\label{lem3.9}
  $\forall \epsilon>0$. Let $g$ be a rapidly decreasing  function  satisfying
 $\sup\limits_{x\in\SR^{n}}\left|x^{\alpha}\partial^{\beta}g\right|<\infty.$
 We denote by $g^{\omega}$ the randomization of $g$ as defined in
(\ref{1.014}).  Then,  there exist $C>0$ and $C_{1}>0$ such that
\begin{eqnarray}
\mathbb{P}\left(\left\{\omega\in \Omega: \left|x^{\alpha}\partial^{\beta}g^{\omega}\right|>M\right\}\right)\leq C_{1}
e^{-\left(\frac{M}{Ce}\right)^{2}}\label{3.034}.
\end{eqnarray}
In particular, take $M=Ce\left[{\rm In}\frac{C_{1}}{\epsilon}\right]^{\frac{1}{2}},$
then, we have
\begin{eqnarray}
\mathbb{P}\left(\left\{\omega\in \Omega: \left|x^{\alpha}\partial^{\beta}g^{\omega}\right|>M\right\}\right)\leq\epsilon\label{3.035}.
\end{eqnarray}
\end{Lemma}
\noindent {\bf Proof}. We firstly show
\begin{eqnarray}
\mathbb{P}\left(\left\{\omega\in \Omega: |x^{\alpha}\partial^{\beta}g^{\omega}|>M\right\}\right)\leq C_{1}
e^{-\left(\frac{M}{Ce}\right)^{2}}\label{3.036}.
\end{eqnarray}
By using Lemmas 3.1, 2.5,  since $g$ is a rapidly decreasing  function which yields
\begin{eqnarray*}
\sum\limits_{|k|\leq2}\int_{|\xi-k|\leq1}\left|\left[(\partial^{\alpha}
\left[\psi(\xi-k))\xi^{\beta}\mathscr{F}_{x}g(\xi)\right]\right]\right|^{2}d\xi\leq C,
\end{eqnarray*}
thus, we have
\begin{eqnarray}
&&\left\|x^{\alpha}\partial^{\beta}g^{\omega}\right\|_{L_{\omega}^{p}(\Omega)}=
\left\|\sum\limits_{k\in Z^{n}}g_{k}(\omega)x^{\alpha}\partial^{\beta}
\psi(D-k)h\right\|_{L_{\omega}^{p}(\Omega)}\nonumber\\
&&=\left\|\sum\limits_{k\in Z^{n}}g_{k}(\omega)\int_{\SR^{n}}e^{ix\xi}\left[-(i\partial^{\alpha})
\left[\psi(\xi-k)(i\xi)^{\beta}\mathscr{F}_{x}g(\xi)\right]\right]\right\|_{L_{\omega}^{p}(\Omega)}\nonumber\\
\nonumber\\&&\leq C\sqrt{p}\sum\limits_{k\in Z^{n}}\left(\int_{\SR^{n}}e^{ix\xi}\left[-(i\partial^{\alpha})
\left[\psi(\xi-k)(i\xi)^{\beta}\mathscr{F}_{x}g(\xi)\right]\right]d\xi\right)^{2}\nonumber\\
&&=C\sqrt{p}\sum\limits_{k\in Z^{n}}\left(\int_{|\xi-k|\leq1}e^{ix\xi}\left[-(i\partial^{\alpha})
\left[\psi(\xi-k)(i\xi)^{\beta}\mathscr{F}_{x}g(\xi)\right]\right]d\xi\right)^{2}\nonumber\\
&&\leq C\sqrt{p}\sum\limits_{k\in Z^{n}}\int_{|\xi-k|\leq1}\left|\left[(\partial^{\alpha}
\left[\psi(\xi-k))(\xi)^{\beta}\mathscr{F}_{x}g(\xi)\right]\right]\right|^{2}d\xi\nonumber\\
&&=C\sqrt{p}\sum\limits_{|k|\leq2}\int_{|\xi-k|\leq1}\left|\left[(\partial^{\alpha}
\left[\psi(\xi-k))\xi^{\beta}\mathscr{F}_{x}g(\xi)\right]\right]\right|^{2}d\xi\nonumber\\&&
\qquad+C\sqrt{p}\sum\limits_{|k|\geq3}\int_{|\xi-k|\leq1}\left|\left[(\partial^{\alpha}
\left[\psi(\xi-k))\xi^{\beta}\mathscr{F}_{x}g(\xi)\right]\right]\right|^{2}d\xi\nonumber\\
&&\leq C\sqrt{p}\sum\limits_{|k|\leq2}\int_{|\xi-k|\leq1}\left|\left[(\partial^{\alpha}
\left[\psi(\xi-k))\xi^{\beta}\mathscr{F}_{x}g(\xi)\right]\right]\right|^{2}d\xi+C\sqrt{p}
\nonumber\\&&\leq C\sqrt{p}.\label{3.037}
\end{eqnarray}
Thus, by Chebyshev inequality and  (\ref{3.037}), we have
\begin{eqnarray}
&&\mathbb{P}\left(\omega\in \Omega:| x^{\alpha}\partial^{\beta}g^{\omega}|>M\right)
\leq \frac{\left\|x^{\alpha}\partial^{\beta}g^{\omega}\right\|_{L_{\omega}^{p}(\Omega)}^{p}}{M^{p}}\leq \frac{(C\sqrt{p})^{p}}{M^{p}}\label{3.038}.
\end{eqnarray}
By using a proof similar to (\ref{3.08}), from (\ref{3.038}), we obtain (\ref{3.034}).

This completes the proof of Lemma 3.9.

\noindent{\bf Remark 4.} From Lemma 3.9, we know  that, if $g$ is a rapidly decreasing  function,
then the randomized function $g^\omega$
is almost surely  a  rapidly decreasing function.
\begin{Lemma}\label{lem3.10}
 $\forall \epsilon>0$ and $\forall \lambda>0$ and  $\|h\|_{L^{2}(\SR^{n})}<\epsilon$
 and we denote by $h^{\omega}$ the randomization of $h$ as defined in
(\ref{1.014}).  Then, there exist $C>0$ and $C_{1}>0$ such that
\begin{eqnarray}
\mathbb{P}\left(\left\{\omega\in \Omega: \|h^{\omega}\|_{L^{2}}>\lambda\right\}\right)\leq C_{1}
e^{-\left(\frac{\lambda}{Ce\|h\|_{L^{2}}}\right)^{2}}\leq C_{1}
e^{-\left(\frac{\lambda}{Ce\epsilon}\right)^{2}}.\label{3.039}
\end{eqnarray}
In particular, take $\lambda=Ce\epsilon\left({\rm In}\frac{C_{1}}{\epsilon}\right)^{\frac{1}{2}},$
and
\begin{eqnarray}
\mathbb{P}\left(\left\{\omega\in \Omega: \|h^{\omega}\|_{L^{2}}>\lambda\right\}\right)\leq C_{1}
e^{-\left(\frac{\lambda}{Ce\|h\|_{L^{2}}}\right)^{2}}\leq\epsilon.\label{3.040}
\end{eqnarray}
\end{Lemma}
\noindent{\bf Proof.} For the proof of (\ref{3.039}), we refer the readers to Lemma 2.2 of \cite{BOP}.
When $\lambda=Ce\epsilon\left({\rm In}\frac{C_{1}}{\epsilon}\right)^{\frac{1}{2}},$ we have
$C_{1}
e^{-\left(\frac{\alpha}{Ce\epsilon}\right)^{2}}=\epsilon.$
We have
\begin{eqnarray}
\lim\limits_{\epsilon\longrightarrow0}\frac{Ce\epsilon\left({\rm In}\frac{C_{1}}{\epsilon}\right)^{\frac{1}{2}}}{\epsilon^{\frac{1}{2}}}=
Ce\lim\limits_{\epsilon\longrightarrow0}\frac{\left({\rm In}\frac{C_{1}}{\epsilon}\right)^{\frac{1}{2}}}{\epsilon^{-\frac{1}{2}}}=0\label{3.038}
\end{eqnarray}

This completes the proof of Lemma 3.10.

\noindent{\bf Remark 5.} From Lemma 3.10, we know that  if $f\in L^{2}(\R^n)$, $n\geq 1$,
then the randomized function $f^\omega$
is almost surely
in $L^{2}(\R^n)$, $n\geq 1$.

\bigskip
\bigskip

\noindent{\large\bf 4. Proof of Theorem 1.1}
\setcounter{equation}{0}
\setcounter{Theorem}{0}

\setcounter{Lemma}{0}

\setcounter{section}{4}
\noindent{\bf Proof of Theorem 1.1}.
When $f\in L^{2}(\R)$, by density theorem which is just Lemma 2.2 in \cite{D}, there exists a rapidly
 decreasing function $g$ such that $f=g+h$, where $\|h\|_{H^{8\epsilon}}<\epsilon.$
We define
\begin{eqnarray}
\Omega_{8}^{c}=\left\{\omega\in \Omega: \left|S_{1}(t)
f^{\omega}-f^{\omega}\right|>\alpha\right\}.\label{4.01}
\end{eqnarray}
Thus, we have
\begin{eqnarray}
\Omega_{8}^{c}\subset \Omega_{9}^{c}\label{4.02}
\cup \Omega_{10}^{c},
\end{eqnarray}
where
\begin{eqnarray}
&&\Omega_{9}^{c}=\left\{\omega\in \Omega: \left|S_{1}(t)
g^{\omega}-g^{\omega}\right|> \frac{\alpha}{2}\right\},\label{4.03}\\
&&\Omega_{10}^{c}=\left\{\omega\in \Omega: \left|S_{1}(t)
h^{\omega}-h^{\omega}\right|> \frac{\alpha}{2}\right\}.\label{4.04}
\end{eqnarray}
Obviously,
\begin{eqnarray}
\Omega_{10}^{c}\subset \Omega_{11}^{c}\cup \Omega_{12}^{c},\label{4.05}
\end{eqnarray}
where
\begin{eqnarray}
&&\Omega_{11}^{c}=\left\{\omega\in \Omega: \left|S_{1}(t)
h^{\omega}\right|> \frac{\alpha}{4}\right\},\label{4.06}\\
&&\Omega_{12}^{c}=\left\{\omega\in \Omega: |h^{\omega}
|> \frac{\alpha}{4}\right\}\label{4.07}.
\end{eqnarray}
From Lemma 3.2, we have
\begin{eqnarray}
\mathbb{P}\left(\Omega_{8}^{c}\right)\leq C_{1}e^{
-\left[\frac{\alpha}{Ce|t|}\right]^{2}}\leq C_{1}e^{
-\left[\frac{\alpha}{Ce|t|}\right]^{2}}\label{4.08}.
\end{eqnarray}
From Lemma 3.3, we have
\begin{eqnarray}
&&\mathbb{P}\left(\Omega_{11}^{c}\right)\leq C_{1}e^{
-\left[\frac{\alpha}{Ce\|h\|_{L^{2}}}\right]^{2}}\leq C_{1}e^{
-\left[\frac{\alpha}{Ce\epsilon}\right]^{2}}\label{4.09}.
\end{eqnarray}
From Lemma 3.8,  we have
\begin{eqnarray}
&&\mathbb{P}\left(\Omega_{12}^{c}\right)\leq C_{1}e^{-\left[\frac{\alpha}{Ce\|h\|_{L^{2}}}\right]^{2}}\leq C_{1}e^{
-\left[\frac{\alpha}{Ce\epsilon}\right]^{2}}\label{4.010}.
\end{eqnarray}
From (\ref{4.08})-(\ref{4.010}),  we have
\begin{eqnarray}
&&\mathbb{P}(\Omega_{8}^{c})\leq \mathbb{P}(\Omega_{9}^{c})+\mathbb{P}(\Omega_{10}^{c})\leq \mathbb{P}(\Omega_{8}^{c})
 +\mathbb{P}( \Omega_{11}^{c})+\mathbb{P}( \Omega_{12}^{c})\nonumber\\&&
 \leq C_{1} e^{-\left[\frac{\alpha}
{C|t|e}\right]^{2}}+2C_{1}e^{
-\left[\frac{\alpha}{Ce\epsilon}\right]^{2}}.\label{4.011}
\end{eqnarray}
When $|t|\leq \epsilon,$    from  (\ref{4.011}), we have
\begin{eqnarray}
\mathbb{P}(\Omega_{8}^{c}) \leq C_{1}e^{
-\frac{\alpha^{2}}{(Ce\epsilon)^{2}}}.\label{4.012}
\end{eqnarray}
Here, $\epsilon$  satisfies $Ce\epsilon\left({\rm In}\frac{C_{1}}{\epsilon}\right)^{\frac{1}{2}}\leq\alpha$.
From (\ref{4.012}), we have
\begin{eqnarray}
\mathbb{P}(\Omega_{8}^{c}) \leq C_{1}e^{-\frac{\alpha^{2}}{(Ce\epsilon)^{2}}}\leq \epsilon.\label{4.013}
\end{eqnarray}
From (\ref{4.013}), we have
\begin{eqnarray}
\mathbb{P}(\Omega_{8}) \geq 1- \epsilon.\label{4.014}
\end{eqnarray}

For the proof of the remainder of Theorem 1.1 can be seen in Lemma 3.11.

This completes the proof of Theorem 1.1.$\hfill\Box$

\bigskip

\bigskip

\noindent {\large\bf 5. Proof of Theorem  1.2}

\setcounter{equation}{0}

 \setcounter{Theorem}{0}

\setcounter{Lemma}{0}

\setcounter{section}{5}
\noindent{\bf Proof of Theorem 1.2}.
When $f\in L^{2}(\R^{n})$, by density theorem which is just Lemma 2.2 in \cite{D},
 there exists a rapidly
 decreasing function $g$ such that $f=g+h$, where $\|h\|_{L^{2}}<\epsilon.$
We define
\begin{eqnarray}
\Omega_{13}^{c}=\left\{\omega\in \Omega: \left|S_{2}(t)
f^{\omega}-f^{\omega}\right|>\alpha\right\}.\label{5.01}
\end{eqnarray}
Thus, we have
\begin{eqnarray}
\Omega_{13}^{c}\subset \Omega_{14}^{c}\label{5.02}
\cup \Omega_{15}^{c},
\end{eqnarray}
where
\begin{eqnarray}
&&\Omega_{14}^{c}=\left\{\omega\in \Omega: \left|S_{2}(t)
g^{\omega}-g^{\omega}\right|> \frac{\alpha}{2}\right\},\label{5.03}\\
&&\Omega_{15}^{c}=\left\{\omega\in \Omega: \left|S_{2}(t)
h^{\omega}-h^{\omega}\right|> \frac{\alpha}{2}\right\}.\label{5.04}
\end{eqnarray}
Obviously,
\begin{eqnarray}
\Omega_{15}^{c}\subset \Omega_{16}^{c}\cup \Omega_{17}^{c},\label{5.05}
\end{eqnarray}
where
\begin{eqnarray}
&&\Omega_{16}^{c}=\left\{\omega\in \Omega: |S_{2}(t)
h^{\omega}|> \frac{\alpha}{4}\right\},\label{5.06}\\
&&\Omega_{17}^{c}=\left\{\omega\in \Omega: |h^{\omega}
|> \frac{\alpha}{4}\right\}\label{5.07}.
\end{eqnarray}
From Lemma 3.4, we have
\begin{eqnarray}
\mathbb{P}\left(\Omega_{13}^{c}\right)\leq C_{1}e^{
-\left[\frac{\alpha}{Ce|t|}\right]^{2}}\leq C_{1}e^{
-\left[\frac{\alpha}{Ce|t|}\right]^{2}}\label{5.08}.
\end{eqnarray}
From Lemma 3.5, we have
\begin{eqnarray}
&&\mathbb{P}\left(\Omega_{16}^{c}\right)\leq C_{1}e^{
-\left[\frac{\alpha}{Ce\|h\|_{L^{2}}}\right]^{2}}\leq C_{1}e^{
-\left[\frac{\alpha}{Ce\epsilon}\right]^{2}}\label{5.09}.
\end{eqnarray}
From Lemma 3.8,  we have
\begin{eqnarray}
&&\mathbb{P}\left(\Omega_{17}^{c}\right)\leq C_{1}e^{-\left[\frac{\alpha}{Ce\|h\|_{L^{2}}}\right]^{2}}\leq C_{1}e^{
-\left[\frac{\alpha}{Ce\epsilon}\right]^{2}}\label{5.010}.
\end{eqnarray}
From (\ref{5.08})-(\ref{5.010}),  we have
\begin{eqnarray}
&&\mathbb{P}(\Omega_{13}^{c})\leq \mathbb{P}(\Omega_{14}^{c})+\mathbb{P}(\Omega_{15}^{c})\leq \mathbb{P}(\Omega_{14}^{c})
 +\mathbb{P}( \Omega_{16}^{c})+\mathbb{P}( \Omega_{17}^{c})\nonumber\\&&
 \leq C_{1} e^{-\left[\frac{\alpha}
{C|t|e}\right]^{2}}+2C_{1}e^{
-\left[\frac{\alpha}{Ce\epsilon}\right]^{2}}.\label{5.011}
\end{eqnarray}
When $|t|\leq \epsilon,$   from  (\ref{5.011}), we have
\begin{eqnarray}
\mathbb{P}(\Omega_{13}^{c}) \leq C_{2}e^{
-\frac{\alpha^{2}}{(Ce\epsilon)^{2}}}.\label{5.012}
\end{eqnarray}
Here, $\epsilon$ satisfies   $Ce\epsilon\left({\rm In}\frac{C_{2}}{\epsilon}\right)^{\frac{1}{2}}\leq\alpha$.
From (\ref{5.012}), we have
\begin{eqnarray}
\mathbb{P}(\Omega_{13}^{c}) \leq C_{2}e^{-\frac{\alpha^{2}}{(Ce\epsilon)^{2}}}\leq \epsilon.\label{5.013}
\end{eqnarray}
From (\ref{5.013}), we have
\begin{eqnarray}
\mathbb{P}(\Omega_{13}) \geq 1- \epsilon.\label{5.014}
\end{eqnarray}
For the proof of the remainder of Theorem 1.2 can be seen in Lemma 3.11.

This completes the proof of Theorem 1.2.$\hfill\Box$

 \bigskip

\bigskip

\bigskip

\bigskip

\noindent {\large\bf 6. Proof of Theorem 1.3}

\setcounter{equation}{0}

 \setcounter{Theorem}{0}

\setcounter{Lemma}{0}

\setcounter{section}{6}

\noindent{\bf Proof of Theorem 1.3}.
When $f\in L^{2}(\R^{n})$, by density theorem which is just Lemma 2.2 in \cite{D}, there exists a rapidly
 decreasing function $g$ such that $f=g+h$, where $\|h\|_{L^{2}(\SR^{n}}<\epsilon.$
We define
\begin{eqnarray}
\Omega_{18}^{c}=\left\{\omega\in \Omega: \left|S_{3}(t)
f^{\omega}-f^{\omega}\right|>\alpha\right\}.\label{6.01}
\end{eqnarray}
Thus, we have
\begin{eqnarray}
\Omega_{18}^{c}\subset \Omega_{19}^{c}\label{6.02}
\cup \Omega_{20}^{c},
\end{eqnarray}
where
\begin{eqnarray}
&&\Omega_{19}^{c}=\left\{\omega\in \Omega: \left|S_{3}(t)
g^{\omega}-g^{\omega}\right|> \frac{\alpha}{2}\right\},\label{6.03}\\
&&\Omega_{20}^{c}=\left\{\omega\in \Omega: \left|S_{3}(t)
h^{\omega}-h^{\omega}\right|> \frac{\alpha}{2}\right\}.\label{6.04}
\end{eqnarray}
Obviously,
\begin{eqnarray}
\Omega_{20}^{c}\subset \Omega_{21}^{c}\cup \Omega_{22}^{c},\label{6.05}
\end{eqnarray}
where
\begin{eqnarray}
&&\Omega_{21}^{c}=\left\{\omega\in \Omega: |S_{3}(t)
h^{\omega}|> \frac{\alpha}{4}\right\},\label{6.06}\\
&&\Omega_{22}^{c}=\left\{\omega\in \Omega: |h^{\omega}
|> \frac{\alpha}{4}\right\}\label{6.07}.
\end{eqnarray}
From Lemma 3.6, we have
\begin{eqnarray}
\mathbb{P}\left(\Omega_{19}^{c}\right)\leq C_{1}e^{
-\left[\frac{\alpha}{Ce|t|}\right]^{2}}\leq C_{1}e^{
-\left[\frac{\alpha}{Ce|t|}\right]^{2}}\label{6.08}.
\end{eqnarray}
From Lemma 3.7, we have
\begin{eqnarray}
&&\mathbb{P}\left(\Omega_{21}^{c}\right)\leq C_{1}e^{
-\left[\frac{\alpha}{Ce\|h\|_{L^{2}}}\right]^{2}}\leq C_{1}e^{
-\left[\frac{\alpha}{Ce\epsilon}\right]^{2}}\label{6.09}.
\end{eqnarray}
From Lemma 3.8,  we have
\begin{eqnarray}
&&\mathbb{P}\left(\Omega_{22}^{c}\right)\leq C_{1}e^{-\left[\frac{\alpha}{Ce\|h\|_{L^{2}}}\right]^{2}}\leq C_{1}e^{
-\left[\frac{\alpha}{Ce\epsilon}\right]^{2}}\label{6.010}.
\end{eqnarray}
From (\ref{4.08})-(\ref{4.010}),  we have
\begin{eqnarray}
&&\mathbb{P}(\Omega_{18}^{c})\leq \mathbb{P}(\Omega_{19}^{c})+\mathbb{P}(\Omega_{20}^{c})\leq \mathbb{P}(\Omega_{19}^{c})
 +\mathbb{P}( \Omega_{21}^{c})+\mathbb{P}( \Omega_{22}^{c})\nonumber\\&&
 \leq C_{1} e^{-\left[\frac{\alpha}
{C|t|e}\right]^{2}}+2C_{1}e^{
-\left[\frac{\alpha}{Ce\epsilon}\right]^{2}}.\label{6.011}
\end{eqnarray}
When $|t|\leq \epsilon,$  from  (\ref{6.011}), we have
\begin{eqnarray}
\mathbb{P}(\Omega_{18}^{c}) \leq C_{2}e^{
-\frac{\alpha^{2}}{(Ce\epsilon)^{2}}}.\label{6.012}
\end{eqnarray}
Here $\epsilon$  satisfies $Ce\epsilon\left({\rm In}\frac{C_{2}}{\epsilon}\right)^{\frac{1}{2}}\leq \alpha$.
From (\ref{6.012}), we have
\begin{eqnarray}
\mathbb{P}(\Omega_{18}^{c}) \leq C_{2}e^{-\frac{\alpha^{2}}{(Ce\epsilon)^{2}}}\leq \epsilon.\label{6.013}
\end{eqnarray}
From (\ref{4.013}), we have
\begin{eqnarray}
\mathbb{P}(\Omega_{18}) \geq 1- \epsilon.\label{6.014}
\end{eqnarray}

For the proof  of the remainder of Theorem 1.3 can be seen in Lemma 3.11.

This completes the proof of Theorem 1.3.$\hfill\Box$

\bigskip
\bigskip

\noindent {\large\bf 7. Proof of Theorem 1.4}

\setcounter{equation}{0}

 \setcounter{Theorem}{0}

\setcounter{Lemma}{0}

\setcounter{section}{7}

In this section, we prove Theorem 1.4.

\noindent{\bf Proof.}For $f\in L^{2}(\R^{n})$, from the density theorem which is just Lemma 2.2 in \cite{D}, we know that there exists
 a decreasing rapidly function $g$   and  $h\in L^{2}(\R^{n})$ with $\|h\|_{L^{2}(\SR^{n})}<\epsilon$ such that $f=g+h.$
Thus, we  have  $f^{\omega}=\sum\limits_{k\in Z^{n}}g_{k}(\omega)P(D-k)f=\sum\limits_{k\in Z^{n}}g_{k}(\omega)P(D-k)(g+h)=
\sum\limits_{k\in Z^{n}}g_{k}(\omega)P(D-k)g+\sum\limits_{k\in Z^{n}}g_{k}(\omega)P(D-k)h=g^{\omega}+h^{\omega}.$
By using a direct computation,  we have
\begin{eqnarray}
&&\mathbb{P}\left(\left\{\omega\in \Omega: \|h^{\omega}\|_{L^{2}}\leq \lambda\right\}\cap \left\{\omega\in
\Omega: \left|x^{\alpha}\partial^{\beta}g^{\omega}\right|\leq M\right\}\right)\nonumber\\
&&= \mathbb{P}\left(\left\{\omega\in \Omega: \|h^{\omega}\|_{L^{2}}\leq \lambda\right\}\right)-\mathbb{P}
\left(\left\{\omega\in \Omega: \|h^{\omega}\|_{L^{2}}\leq \lambda\right\}\cap \left\{\omega\in \Omega:
 \left|x^{\alpha}\partial^{\beta}g^{\omega}\right|> M\right\}\right)\nonumber\\
&&\geq \mathbb{P}\left(\left\{\omega\in \Omega: \|h^{\omega}\|_{L^{2}}\leq \lambda\right\}\right)-\mathbb{P}
\left(\left\{\omega\in \Omega: \left|x^{\alpha}\partial^{\beta}g^{\omega}\right|> M\right\}\right)\nonumber\\&&\geq1-\epsilon-\epsilon=1-2\epsilon.
\end{eqnarray}

This completes the proof of Theorem 1.4.

\bigskip

\leftline{\large \bf Acknowledgments}
Wei Yan was supported by NSFC grants (No. 11771127) and the Young core Teachers program of
 Henan province under grant number 5201019430009,  Jinqiao Duan was supported by the NSF grant (No. 1620449),
and NSFC grants (No. 11531006, No. 11771449), Yongsheng Li supported by NSFC grants (No. 11571118)
 and   Meihua Yang was supported by
NSFC grants (No. 11571125).

\end{document}